\documentclass[12pt]{article}

\usepackage{amsmath,epsfig,epsf,psfrag,natbib,bbm}
\usepackage{amssymb}
\usepackage[latin1]{inputenc}
\usepackage{rotating}
\usepackage{amsbsy}
\usepackage{lscape}
\usepackage{color,latexsym,amsfonts,amssymb,amsthm,amscd,amsmath,graphicx}

\newcommand{\be}{\begin{equation}}
\newcommand{\ee}{\end{equation}}
\newtheorem{theorem}{Theorem}
\newtheorem{lemma}{Lemma}

\newcommand{\var}{\mbox{\sc var}}
\newcommand{\E}{\mathbb E}
\newcommand{\cov}{{\mbox{\sc cov}}}

\newcommand{\btheta}{\mbox{\boldmath $\theta$}}

\renewcommand{\ss}{\mbox{\bf s}}
\renewcommand{\tt}{\mbox{\bf t}}

\newcommand{\bx}{\mbox{\bf x}}

\newcommand{\bw}{\mbox{\bf w}}

\newcommand{\ind}{\mbox{$\mathbbm{1}$}}

\def\conD{\stackrel{\cal D} \to}

\begin{document}
\title{Test for homogeneity with unordered paired observations}

\author{Jiahua Chen, Pengfei Li, Jing Qin, and  Tao Yu}

\date{}
\maketitle

\vspace{-5mm}

\begin{abstract}
In some applications, an experimental unit is composed of
two distinct but related subunits. The response from such a unit is $(X_{1}, X_{2})$
but we observe only $Y_1 = \min\{X_{1},X_{2}\}$ and $Y_2 = \max\{X_{1},X_{2}\}$, i.e.,
the subunit identities are not observed.
We call $(Y_1, Y_2)$ unordered paired observations. Based on unordered paired observations $\{(Y_{1i}, Y_{2i})\}_{i=1}^n$,
we are interested in whether the
marginal distributions for $X_1$ and $X_2$ are identical.
Testing methods are available in the literature under the assumptions that
$\var(X_1) = \var(X_2)$ and $\cov(X_1, X_2) = 0$.
However, by extensive simulation studies, we observe that when one or both assumptions are violated,
these methods have inflated type I errors or much lower
powers.
In this paper, we study the likelihood ratio test statistics
for various scenarios and explore their limiting distributions
without these restrictive assumptions.
Furthermore, we develop Bartlett correction formulae for these
statistics to enhance their precision when
the sample size is not large.
Simulation studies and real-data examples are used to illustrate
the efficacy of the proposed methods.
\end{abstract}

\section{Introduction}
In some applications, an experimental unit is made of
two distinct but related subunits. The response from such a unit is $(X_{1}, X_{2})$
but we observe only $Y_1 = \min\{X_{1},X_{2}\}$ and $Y_2 = \max\{X_{1},X_{2}\}$; that is,
the subunit identities are not observed or unobservable.
We call $(Y_1, Y_2)$ unordered paired observations.
We assume that $(X_{1i}, X_{2i})^\tau$, for $i=1, \ldots, n$,
are independent and identically distributed (i.i.d.)
normal random vectors:
\begin{equation}
\label{bn}
\left(
\begin{array}{c}
X_{1i}\\
X_{2i}
\end{array}
\right)
\sim
N\left(
\left(
\begin{array}{c}
\mu_1\\
\mu_2\\
\end{array}
\right),
\left(
\begin{array}{cc}
\sigma_1^2&\rho\sigma_1\sigma_2\\
\rho\sigma_1\sigma_2&\sigma_2^2\\
\end{array}
\right)
\right).
\end{equation}
We say that $\{(Y_{1i}, Y_{2i})\}_{i=1}^n$ are {\it uncorrelated} when $\rho =0$
and {\it correlated} when $\rho \neq 0$.
This paper studies the homogeneity testing of the marginal distributions of $X_{1i}$ and $X_{2i}$:
\begin{eqnarray}
H_0: (\mu_1,\sigma_1^2) = (\mu_2,\sigma_2^2)
\quad \mbox{versus}
\quad
H_a: (\mu_1,\sigma_1^2) \neq (\mu_2,\sigma_2^2).
\label{hypothesis}
\end{eqnarray}

Unordered paired data occur in many applications,
and there is a long research history. For instance,
\cite{hinkley1973} analyzed such a data set from human genetics.
The genetic blueprint of an individual is contained in 23 pairs of chromosomes.
Each member of the pair is inherited from the corresponding
chromosome pair of a parent.
If we do not know the chromosome correspondences
between the offspring and the parents,
we lose the parental identities and end up with unordered paired observations.
\cite{olkin1995} provide more examples.
In visual acuity studies, we may record only
a subject's extreme acuities (the ``best" and ``worst" acuities) without
recording the corresponding eyes.
In twin experiments, we obtain unordered paired observations
without a label for each member of a twin pair;
see \cite{ernst1996} and \cite{shekar2006} and the references therein.
Furthermore, unordered data of a higher dimension are 
collected in various scientific disciplines.
For example, \cite{davies1988} provided an example of unordered data
of dimension $k$.
In the interim analysis of a double-blinded clinical trial of $k$ treatments,
we get the $k$ order statistics without knowledge of the corresponding treatments;
see also \cite{van2005} and \cite{miller2009}.
In diffusion tensor (DT) brain imaging (see \cite{yu2013} and the references therein),
the eigenvalues of the DT estimates for each brain voxel are viewed as
unordered triples.

With unordered paired observations, a fundamental question is
whether or not $X_{1i}$ and $X_{2i}$ have the same distribution.
Under Model (\ref{bn}), this is equivalent to testing the hypothesis specified
in (\ref{hypothesis}).
\cite{hinkley1973} proposed a likelihood ratio test (LRT) procedure
under the assumption that $\rho = 0$ and $\sigma_1^2=\sigma_2^2$.
\cite{li2011} investigated this problem in a semiparametric setup.
Other approaches can be found in
\cite{mooreII1973}, \cite{lauder1977}, \cite{mooreII1979},
\cite{carothers1981}, \cite{efron1971}, and \cite{qin2005}, among others.
All these works assume that $X_{1i}$ and $X_{2i}$
are independent with equal variance.
These assumptions may not hold in applications,
and they can be severely violated, as evidenced by the
examples in Section \ref{section-real-data}.
Ignoring the dependence structure and/or imposing an incorrect equal-variance assumption
can lead to unreliable inference conclusions:
the type I error may be severely inflated or the power
markedly decreased.

This paper focuses on tests for (\ref{hypothesis}).
In particular, we study the LRT in four scenarios:
(1) $\rho = 0$ and $\sigma_1^2 = \sigma_2^2$;
(2) $\rho = 0$;
(3) $\sigma_1^2 = \sigma_2^2$;
and
(4) no assumption on $\rho$, $\sigma_1^2$, and $\sigma_2^2$.

Investigating the asymptotic behavior  of  these LRT statistics is  technically challenging.
The well-developed theory
\citep{wilks1938, chernoff1954, self1987, drton2009} is not applicable 
because of the undesirable mathematical properties
(see (\ref{prop1}) in Section \ref{section-2})
of the log-likelihood function.
In addition, an important byproduct of the theory for the corresponding
LRT statistics is the asymptotic behavior of the maximum likelihood estimators (MLEs)
for $(\mu_1, \mu_2, \sigma^2_1,\sigma_2^2)$.
Interestingly,  we have shown that the asymptotic behavior depends on whether
$\rho=0$ is known or $\rho$ is unknown.
The convergence rates of these parameter estimates
depend on the scenario.

We observe that the limiting distributions of the LRT statistics under $H_0$ are not
sufficiently accurate approximations to their finite-sample distributions
when $n$ is not large.
To enhance the approximation precision of the limiting distributions,
we adjust the statistics based on the Bartlett correction
\citep{bartlett1937,lawley1956}. Simulation results confirm
the efficacy of the adjustment.

We organize the rest of the paper as follows.
Section \ref{section-2} introduces the LRT statistics for (\ref{hypothesis})
and studies their asymptotic behavior under $H_0$.
Section \ref{section-3} presents the adjusted limiting distributions
of our statistics for data of limited sample size.
Section \ref{section-4} contains simulation studies, and
Section \ref{section-real-data} gives real-data examples.
The technical details are relegated to Section \ref{appendix}.

\section{Main Results\label{section-2}}
The LRT is an essential tool in statistical inference,
especially under the parametric model assumption;
see \cite{wilks1938, chernoff1954, self1987, drton2009},
and the references therein.
In this section, we present LRT statistics
and study their properties for testing (\ref{hypothesis})
under model assumptions on
$\rho$ and whether or not $\sigma_1^2 =\sigma_2^2$.

We first derive the log-likelihood function with unordered paired observations.
{For any $y_1 < y_2$,  we have
\begin{eqnarray*}
P(Y_1 \leq y_1, Y_2 \leq y_2)
&=& P\left(\{X_1 \leq y_1, X_2 \leq y_2\} \cup  \{X_1 \leq y_2, X_2 \leq y_1\}\right)\\
&=&P(X_1 \leq y_1, X_2 \leq y_2 ) + P( X_1 \leq y_2, X_2 \leq y_1 ) \\
&& -P(\{X_1 \leq y_1, X_2 \leq y_2\}\cap \{X_1 \leq y_2, X_2 \leq y_1\})\\
&=& P(X_1 \leq y_1, X_2 \leq y_2 ) + P( X_1 \leq y_2, X_2 \leq y_1 ) -P(X_1 \leq y_1, X_2 \leq y_1).
\end{eqnarray*}
}
Therefore, the joint density function of $(Y_1, Y_2)$ is given by
\[
\phi (y_{1}, y_{2}; \btheta)
+
\phi (y_{2}, y_{1}; \btheta),
\]
where $\phi(x_1, x_2; \btheta)$ denotes the bivariate normal density function
with parameters $\btheta = (\mu_1,\mu_2,\sigma_1,\sigma_2,\rho)^\tau$
specified in \eqref{bn}.
The log-likelihood function based on $\{(Y_{1i}, Y_{2i})\}_{i=1}^n$ and
Model \eqref{bn} is:
\begin{equation}
\label{loglik1}
\ell_n(\btheta)
=
\sum_{i=1}^{n}
\log\{ \phi (Y_{1i}, Y_{2i}; \btheta)
+
\phi (Y_{2i}, Y_{1i}; \btheta) \}.
\end{equation}
This likelihood function is the basis for our subsequent development.

\subsection{Unordered uncorrelated paired data}
\label{sec-2-1}

In this section, we assume that $\rho = 0$ is known; problem (\ref{hypothesis}) is reduced to
$H_0: \mu_1 = \mu_2, \sigma_1= \sigma_2$.
We define
\begin{eqnarray*}
\hat \btheta
&=& \arg\sup_{\btheta}\{ \ell_{n}(\btheta):  \rho=0\},
\\
\tilde \btheta
&=& \arg\sup_{\btheta}\{ \ell_{n}(\btheta): \sigma_1= \sigma_2=\sigma, \rho=0\},
\\
\check {\btheta}
&=& \arg\sup_{\btheta}\{ \ell_{n}(\btheta): (\mu_1, \sigma_1) = (\mu_2, \sigma_2), \rho=0\},
\end{eqnarray*}
and we use the notational convention that the entries of $\hat \btheta$ are
$\hat \mu_1$, $\hat \mu_2$, and so on.
Note that  $\hat \btheta$, $\tilde \btheta$, and
$\check \btheta$ are MLEs of $\btheta$ under various
constraints.
The LRT statistics for testing the null hypothesis  \eqref{hypothesis}
against two alternatives, specified by
$\sigma_1 = \sigma_2$ and $\sigma_1 \neq \sigma_2$
respectively, are given by
\be
\label{Rn1}
R_{n,1}
=
2 \{ \ell_{n}(\tilde \btheta) - \ell_{n}(\check \btheta)\}, ~~~
R_{n,2}
=
2 \{ \ell_{n}(\hat \btheta) - \ell_{n}(\check \btheta)\}.
\ee

Theorem \ref{theorem-1} below establishes the asymptotic distributions
of $R_{n,1}$ and $R_{n,2}$ as well as the convergence rates of
$\tilde \btheta$ and $\hat \btheta$ under $H_0$.
For presentational continuity, we relegate its proof to Section \ref{appendix}.
Let $\conD$ denote ``convergence in distribution."
We use $0.5\chi^2_0+0.5\chi^2_1$ for  an equal mixture of $\chi_0^2$ and $\chi_1^2$,
with $\chi_0^2$ being the distribution with a point mass at zero.

\begin{theorem}
\label{theorem-1}
Assume Model (\ref{bn}) and $\rho = 0$.
Under $H_0$, as $n\to \infty$, we have
\begin{enumerate}
\item[(a)]
$(\tilde \mu_1-\mu_0)^2, (\tilde \mu_2-\mu_0)^2$, and
$\tilde \sigma-\sigma_0$ are all of order $O_p(n^{-1/2})$,
and
\[
R_{n,1}\conD0.5\chi^2_0+0.5\chi^2_1;
\]
\item[(b)]
$(\hat \mu_j-\mu_0)^2$, $(\hat\sigma_j-\sigma_0)^2$
for $j=1, 2$ are all of order $O_p(n^{-1/2})$, and
\[
R_{n,2}\conD R\equiv\sup_{x_1,x_2} \left\{ 2\bx^\tau \bw-\bx^\tau \bx\right\},
\]
where
$\bx^\tau=(x_1^2,x_2^2,2x_1x_2)$ and
$\bw^\tau=(w_1,w_2, w_3)$ with
$w_1, w_2, w_3$ being three {i.i.d.}  $N(0,1)$ random variables.
\end{enumerate}
\end{theorem}

Deriving the asymptotic null distributions of $R_{n,1}$ and $R_{n,2}$
is technically challenging.
We make the following comments.
Let $\mu=(\mu_1+\mu_2)/2$ and $\Delta=(\mu_1-\mu_2)/2$
so that $\mu_1=\mu+\Delta$ and $\mu_2=\mu-\Delta$; we have
\begin{eqnarray}
\label{prop1}
\frac{\partial \ell_n(\mu+\Delta, \mu - \Delta, \sigma_1, \sigma_2,\rho)}
{\partial \Delta} \Big |_{\Delta= 0, \sigma_1 = \sigma_2} =0.
\end{eqnarray}
This fact implies that the Fisher information matrix of
$\btheta $
under the null hypothesis degenerates and undermines
the basis for the elegant classical results
\citep{wilks1938, chernoff1954, self1987, drton2009}.
The crucial step in obtaining the asymptotic null distribution
of the LRT is a quadratic approximation in $\hat \btheta - \btheta$
to the log-likelihood ratio function.
Following this path, we need to consider a fourth-order Taylor expansion to obtain a quadratic approximation in
$(\hat \btheta - \btheta)^2$ and so on.
Fortunately, we find that the sandwich technique of \cite{chen2001}
and \cite{chenjrss2001}
overcomes the technical obstacles caused by (\ref{prop1}).

\subsection{Unordered correlated pair data}

In this section, we study the LRTs for  (\ref{hypothesis}) with $\rho$ being an unknown parameter.
Define
\begin{eqnarray*}
\hat \btheta^*
&=& \arg\sup_{\btheta}\{ \ell_{n}(\btheta)\}, \\
\tilde \btheta^*
&=& \arg\sup_{\btheta}\{ \ell_{n}(\btheta): \sigma_1= \sigma_2=\sigma\}, \\
\check \btheta^*
&=& \arg\sup_{\btheta}\{ \ell_{n}(\btheta): (\mu_1, \sigma_1) = (\mu_2, \sigma_2)\}.
\end{eqnarray*}
Similarly to the strategy for (\ref{Rn1}), we define the LRT statistics for  (\ref{hypothesis}) with $\rho$ being an unknown parameter:
\begin{eqnarray*}
R^*_{n,1}
=
2 \{ \ell_{n}(\tilde \btheta^*)  - \ell_{n}(\check \btheta^*) \}, ~~~
R^*_{n,2}
=
2 \{ \ell_{n}(\hat \btheta^*)  - \ell_{n}(\check \btheta^*) \}.
\end{eqnarray*}

Theorem \ref{theorem-2} below establishes the asymptotic distributions of
$R_{n,1}^*$ and $R_{n,2}^*$ as well as the convergence rates of
$\tilde \btheta^*$ and $\hat \btheta^*$ under
their respective $H_0$.
The proof is given in Section \ref{appendix}.

\begin{theorem} \label{theorem-2}
Assume Model (\ref{bn}) but do not assume $\rho = 0$.
Under $H_0$, as $n\to \infty$,  we have
\begin{enumerate}
\item[(a)]
$(\tilde \mu_1^*-\mu_0)^2$,
$(\tilde \mu_2^*-\mu_0)^2$,
$(\tilde\sigma^*-\sigma_0)$, and
$(\tilde\rho^*-\rho_0)$
are all of order  $O_p\left(n^{-1/4}\right)$, and
$$
R_{n,1}^*\conD0.5\chi^2_0+0.5\chi^2_1;
$$
\item[(b)]
$(\hat \mu_1^*-\mu_0)^2$,
$(\hat \mu_2^*-\mu_0)^2$,
$\hat\sigma^*_1-\sigma_0$,
$\hat\sigma^*_2-\sigma_0$, and
$\hat\rho^*-\rho_0$
are all of order $O_p\left(n^{-1/4}\right)$, and
$$
R_{n,2}^*\conD R^*\equiv \max\{w_1^2+(w_2^+)^2, w_1^2+(w_3^+)^2\},
$$
where $w_1$, $w_2$, and $w_3$ are three {i.i.d.}  $N(0,1)$ random variables.
\end{enumerate}
\end{theorem}

The limiting cumulative distribution function (c.d.f.)
of $R^*_{n,2}$ is given by:
\[
P(R^*\leq x) = P \Big( \max\{w_1^2+(w_2^+)^2, w_1^2+(w_3^+)^2\} \leq x \Big )
=
\int_{0}^x \Phi^2(\sqrt{x-y})(2\pi y)^{-1/2}\exp(-y/2) dy
\]
for $x\geq 0$ with $\Phi(\cdot)$ being the c.d.f. of the standard normal distribution.
We use this expression to evaluate
the asymptotic quantile and the p-value for the corresponding test.

\section{Adjusted Limiting Distributions}
\label{section-3}

One drawback of the general asymptotic results is that they may
offer poor approximations to the corresponding finite-sample distributions. 
The convergence rates of the parameter estimators
given in Theorems \ref{theorem-1} and \ref{theorem-2}
 are much lower than those of the MLEs from the regular parametric models.
This adversely affects the approximation accuracy of the asymptotic distributions to
 the finite-sample distributions of the LRT statistics.
To improve the approximation precision when $n$ is not very large,
we use the Bartlett correction.
Suppose the limiting distribution of a statistic $T_n$ is given by $F(x)$.
We may search for  a sequence of c.d.f.s $F_n(x) \to F(x)$ such that $F_n(x)$
and $T_n$ have the same first moment up to order $O\left(n^{-1}\right)$.
This idea was pioneered by \cite{bartlett1937} and generalized by \cite{lawley1956}.

In this spirit, we search for  accurate approximate distributions
for $R_{n,1}$, $R_{n,2}$, $R_{n,1}^*$, and $R_{n,2}^*$ as follows.
Recall that $R$ and $R^*$ are the limiting distributions of $R_{n,2}$ and
$R^*_{n,2}$. 
Let
\begin{eqnarray*}
F_{n1} & = (1-p_n)\chi^2_0+p_n\chi^2_1, 
\label{R_n1-adj}\\
F_{n2} &= r_n R, \quad
 \label{R_n2-adj}\\
F_{n1}^* &= (1-p_n^*)\chi^2_0+p_n^*\chi^2_1,
\label{R_n1-star-adj}\\
F_{n2}^* & =  r_n^* R^*
\label{R_n2-star-adj}.
\end{eqnarray*}
We need to find $p_n$, $r_n$, $p^*_n$, and $r^*_n$
so that the above distributions have first moments very close to the first moments
of their corresponding test statistics for a wide range of $n$ values.
High-order asymptotic techniques can be used,
but they may involve complicated analytical tools with
little assurance of the quality of the end products.
The computer experiment approach of \cite{chen2011tuning}
is more effective and practical, and it matches the spirit of the data science.

The experiment works as follows. 
{We consider a sufficiently wide range of values for $n$.
For each $n$, 
 we simulate a large number of data sets, with each data set composed of $n$  i.i.d. unordered paired observations. 
 Due to the invariance property of the LRT statistics, 
each data set is generated from the standard bivariate normal distribution.
Based on these data sets, we obtain the simulated first moments of $R_{n,1}$, $R_{n,2}$, $R_{n,1}^*$,
and $R_{n,2}^*$.
We choose $p_n$ so that   the simulated first moment of $R_{n,1}$ matches the first moment of $F_{n1}$. 
We then look for a regression model for $p_n$ versus $n$. 
Similar procedures are applied to obtain regression models for $r_n$, $p^*_n$, and $r^*_n$.

Specifically, let us take $R_{n,1}$ for ease of illustration:
\begin{enumerate}
\item[]Step 1.
For every $n$ in  $\{10, 20, \ldots, 100\}$, generate $N= 50,000$ data sets of size $n$.
\item[]Step 2.
Obtain $N$ values of $R_{n,1}$  and therefore  its simulated first moment, denoted $\hat p_n$.
Match $\hat p_n$ with the first moment of $F_{n1}$ to find $p_n=\hat p_n$. 

\item[]Step 3. Fit a regression model to $(n,p_n)$ with $p_n $ being the response and $n$ being the covariate.
\end{enumerate}

We postulate the following nonlinear but parametric regression models:
\begin{eqnarray}
p_n &=&
 	0.5+a n^{-b}+\epsilon_n \label{p_n_adj} \\
r_n &=&
	1+an^{-b}+\epsilon_n \label{r_n_adj} \\
p_n^* &=&
 	0.5+an^{-b}+\epsilon_n \label{p_n_star_adj}\\
r_n^* &=&
	1+an^{-b}+\epsilon_n, \label{r_n_star_adj}
 \end{eqnarray}
with $a$ and $b$ being regression parameters, and $\epsilon_n$ accounting for imperfect fit.
Applying Steps 1--2 outlined above leads to the
$ p_n$, $ r_n$,  $  p_n^*$, and $ r_n^*$ values
 in Table \ref{pn}.
Fitting the nonlinear regression models (\ref{p_n_adj})--(\ref{r_n_star_adj})
to the data in Table \ref{pn}  gives us the fitted values of $a$ and $b$. 
With these values, we calculate the approximate p-values with the following adjusted limiting distributions:
\begin{eqnarray*}
(0.5-1.440n^{-0.676})\chi^2_0+(0.5+1.440n^{-0.676}) \chi^2_1 \quad
&\mbox{for}& R_{n,1}, \label{R_n1-adj2}\\
 (1+4.589n^{-1.163} )R  \quad &\mbox{for}& R_{n,2}, \label{R_n2-adj2}\\
(0.5-1.332n^{-0.492})\chi^2_0+(0.5+1.332n^{-0.492}) \chi^2_1 \quad
&\mbox{for}& R_{n,1}^*, \label{R_n1-star-adj2}\\
 (1+6.325n^{-1.176} )R^*\quad &\mbox{for}& R_{n,2}^*. \label{R_n2-star-adj2}
\end{eqnarray*}
We have implemented the four LRT statistics with the proposed adjusting limiting distributions
in an \texttt{R} package; it is available upon request. 
}

\begin{table}[http]
\caption{Values of $p_n$, $ r_n$,  $ p_n^*$, and $ r_n^*$ via computer
experiments}
\label{pn}
\renewcommand{\arraystretch}{.8}
\begin{center}
\tabcolsep=2mm
\begin{tabular}{lcccccccccc}
\hline
$n$&     10   &  20   &  30 &  40  &  50 &  60  &  70 &  80  &  90  &  100 \\ \hline
$p_n$&0.809&0.681&0.634&0.627&0.596&0.587&0.585&0.587&0.568&0.568\\
$r_n$&1.312&1.150&1.092&1.070&1.046&1.028&1.030&1.032&1.016&1.012\\
$p_n^*$&0.932&0.801&0.749&0.721&0.687&0.674&0.669&0.651&0.649&0.645\\
$r_n^*$&1.417&1.194&1.129&1.090&1.062&1.040&1.038&1.028&1.022&1.018\\
\hline
\end{tabular}
\end{center}
\end{table}


\section{Simulation Studies}\label{section-4}

\subsection{Data generation} \label{section-4-1}

Because of the invariance property, we need only study
the LRT tests based on data generated from distributions
with standardized parameter values.

To examine the sizes of the tests, we simulate at $\mu_1 = \mu_2 = 0$
and $\sigma_1 = \sigma_2 = 1$ in \eqref{bn}.
We study five cases corresponding to $\rho = -0.5, -0.25, 0, 0.25$, and $0.5$.
To compare the powers of the tests, we set $\mu_1 = 0$, $\sigma_1 = 1$,
and form 20 cases as combinations of $\mu_2 = 1.0, 1.5$,
$\sigma_2 = 1.0, 0.5$ and $\rho = -0.5, -0.25, 0, 0.25, 0.5$.

In each case, we generate $(X_1, X_2)$ from model \eqref{bn}
with one of the above parameter settings.
Then, we obtain $Y_1 = \min \{X_1, X_2\}$ and $Y_2 = \max\{X_1,X_2\}$.
We repeat the process to obtain $n$ unordered pairs $(Y_1, Y_2)$.

Based on each set of $n$ unordered pairs,
we compute the values of $R_{n,1}$, $R_{n,2}$, $R_{n,1}^*$,
and $R_{n,2}^*$ and carry out the tests for $H_0$
without checking that the model for generating the data satisfies
the conditions for the tests.
We record the rejection rates based on $50,000$
repetitions; the results are presented
in the next section.

\subsection{Results}
\label{section-4-2}
We calculate the rejection rate of each test
at the significance levels $\alpha = 10\%, 5\%$, and $1\%$.
The rejection percentages under the null models
are summarized in Table \ref{type1}.
\begin{table}[!ht]
 \caption{Simulated Type I errors (\%) of LRTs based on limiting distributions/adjusted limiting distributions}
 \label{type1}
 \renewcommand{\arraystretch}{.8}
\begin{center}
\tabcolsep=2mm
\begin{tabular}{r| rrr|rrr}
\hline
Levels&  10\% & 5\% & 1\% & 10\% & 5\% & 1\% \\
\hline
&\multicolumn{3}{|c|}{$n=25$}&\multicolumn{3}{c}{$n=75$}\\ \hline
&\multicolumn{6}{c}{$\rho=0$}\\
$R_{n,1}$   &13.7/10.7&7.3/5.7&1.8/1.4&11.3/\ 9.9&5.9/5.1&1.3/1.2\\
$R_{n,2}$   &12.9/10.6&6.9/5.2&1.6/1.0&10.8/10.2&5.6/5.2&1.2/1.0\\
$R_{n,1}^*$ &15.9/10.5&8.1/5.5&1.8/1.1&13.4/10.4&7.0/5.5&1.5/1.1\\
$R_{n,2}^*$ &13.5/10.1&7.4/5.0&1.8/1.1&11.1/10.1&5.9/5.2&1.2/1.0\\
\hline
&\multicolumn{6}{c}{$\rho=0.25$}\\
$R_{n,1}$   &1.2/0.8&0.5/0.3&0.1/0.1&0.1/0.0&0.0/0.0&0.0/0.0\\
$R_{n,2}$   &3.8/3.0&1.9/1.4&0.4/0.3&1.8/1.7&0.7/0.7&0.1/0.1\\
$R_{n,1}^*$ &15.9/10.5&8.1/5.5&1.8/1.1&13.4/10.4&7.0/5.5&1.5/1.1\\
$R_{n,2}^*$ &13.5/10.1&7.4/5.0&1.8/1.1&11.1/10.1&5.9/5.2&1.2/1.0\\
\hline
&\multicolumn{6}{c}{$\rho=0.5$}\\
$R_{n,1}$   &0.0/0.0&0.0/0.0&0.0/0.0&0.0/0.0&0.0/0.0&0.0/0.0\\
$R_{n,2}$   &0.7/0.5&0.3/0.2&0.0/0.0&0.1/0.1&0.0/0.0&0.0/0.0\\
$R_{n,1}^*$ &15.9/10.5&8.1/5.5&1.8/1.1&13.4/10.4&7.0/5.5&1.5/1.1\\
$R_{n,2}^*$ &13.5/10.1&7.4/5.0&1.8/1.1&11.1/10.1&5.9/5.2&1.2/1.0\\
\hline
&\multicolumn{6}{c}{$\rho=-0.25$}\\
$R_{n,1}$   &53.7/47.2&38.6/33.0&15.2/12.7&83.1/80.9&71.6/69.1&43.6/41.3\\
$R_{n,2}$   &39.0/34.0&25.5/21.2&8.6/6.2&67.6/66.3&53.6/52.0&27.3/25.6\\
$R_{n,1}^*$ &15.9/10.5&8.1/5.5&1.8/1.1&13.4/10.4&7.0/5.5&1.5/1.1\\
$R_{n,2}^*$ &13.5/10.1&7.4/5.0&1.8/1.1&11.1/10.1&5.9/5.2&1.2/1.0\\
\hline
&\multicolumn{6}{c}{$\rho=-0.5$}\\
$R_{n,1}$   &92.6/89.9&84.5/80.5&57.5/52.4&100.0/99.9&99.9/99.8&98.5/98.3\\
$R_{n,2}$   &80.1/76.2&67.1/61.3&37.3/30.2&99.7/99.6&99.0/98.9&94.5/93.9\\
$R_{n,1}^*$ &15.9/10.5&8.1/5.5&1.8/1.1&13.4/10.4&7.0/5.5&1.5/1.1\\
$R_{n,2}^*$ &13.5/10.1&7.4/5.0&1.8/1.1&11.1/10.1&5.9/5.2&1.2/1.0\\
\hline
\end{tabular}
\end{center}
\end{table}

When $\rho=0$, $X_{1}$ and $X_{2}$ are simulated to be independent.
The assumptions for all the LRTs, $R_{n,1}$, $R_{n,2}$, $R_{n,1}^*$, and $R_{n,2}^*$,
are satisfied.
However, as shown in the first section of Table \ref{type1},
if their limiting distributions are applied without adjustment,
the resulting tests are inaccurate:
their type I errors markedly exceed the nominal significance levels.
The adjustment proposed in Section \ref{section-3} is very helpful.
After the adjustment,
the type I errors of all the tests are close to the nominal levels.
The precision is impressive since the adjustment works well even when $n$ is
as small as $25$.

When $\rho = \pm 0.25$ or $\pm 0.5$, the model assumptions for $R_{n,1}$ and $R_{n,2}$
are violated. When we apply the tests, the type I errors are either near zero
when $\rho = 0.25$ or $0.5$ or seriously inflated when $\rho = -0.25$ or $-0.5$.
In contrast, because of their invariance property, $R_{n,1}^*$ and $R_{n,2}^*$ continue to perform well:  with their limiting distributions adjusted, they have
satisfactory precision in the type I errors.

\begin{figure}[!ht]
\caption{
Simulated type I errors (\%) at the 5\% significance level  when $100 \leq n \leq 1500$ and $\rho=0$.
The solid and dashed lines are the rates before and after the adjustments, respectively.
}
\centerline{ \includegraphics[scale=0.6]{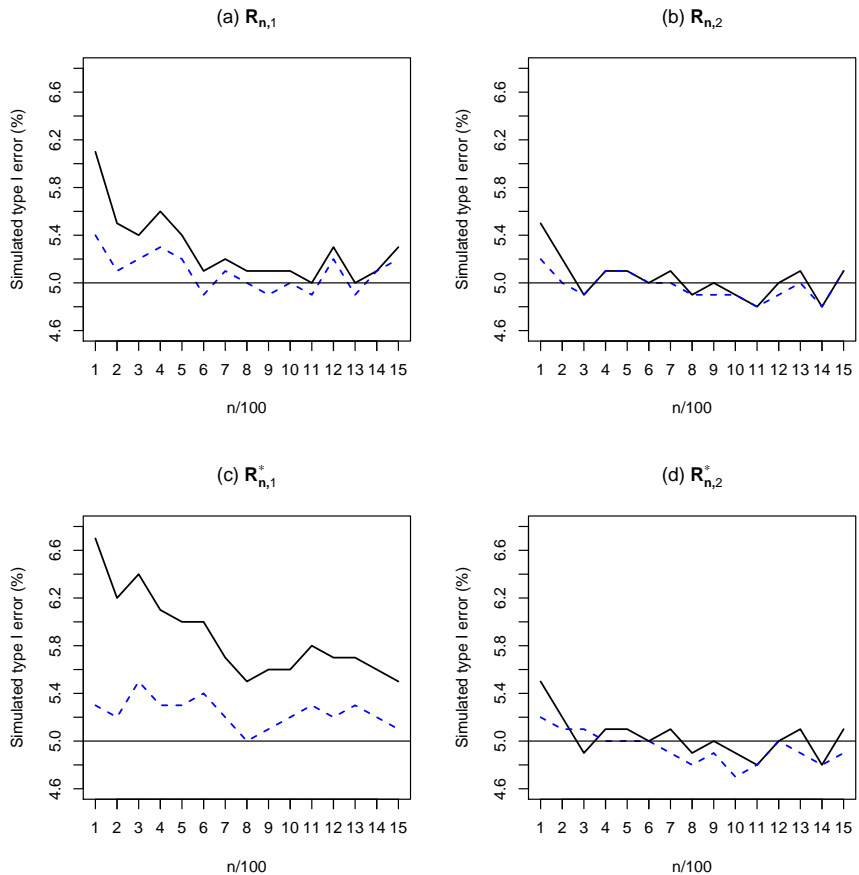}}
\label{typeI_5}
\end{figure}

To further illustrate the effects of the adjustment on the limiting distributions,
Figure \ref{typeI_5} presents the type I errors (\%) of
our LRTs at the 5\% significance level when $100 \leq n \leq1500$ and $\rho=0$.
The trends for the 10\% and 1\% significance levels  are similar and are omitted.
The plots show that the type I errors of $R_{n,1}$, $R_{n,2}$
after the adjustment are within a $0.2\%$ band of the nominal level for large $n$ and a $0.4\%$ band otherwise; similar results are observed for $R_{n,1}^*$. For $R_{n,1}^*$,  the approximation accuracy shows no clear improvement as $n$ increases, but the type I errors are between  5\% and 5.4\%, which is sufficiently accurate for typical applications.

Next, we compare the powers of $R_{n,1}$, $R_{n,2}$, $R_{n,1}^*$, and $R_{n,2}^*$
under the alternatives. All combinations of $n$, $\rho$, $\mu$, and $\sigma$ are
incorporated, as described in Section \ref{section-4-1}.
Their powers, summarized in Table \ref{power},
are computed at the 5\% significance level based on the
adjusted limiting distributions.
We observe that when $\rho=0$, $R_{n,1}$ and  $R_{n,2}$ have higher powers than
$R_{n,1}^*$ and $R_{n,2}^*$;
when $\rho=0.25$, $R_{n,1}$ and  $R_{n,2}$ have higher powers in most cases;
when $\rho$ is increased to 0.5, $R_{n,1}^*$ and $R_{n,2}^*$ are much more powerful;
when $\rho=-0.25$ and $-0.5$, $R_{n,1}$ and  $R_{n,2}$ are more powerful,
but at the cost of the inflated type I errors reported in Table \ref{type1};
a test with a markedly inflated type I error is generally not recommended.

\begin{table}[!ht]
 \caption{Powers (\%) of  $R_{n,1}$, $R_{n,2}$, $R_{n,1}^*$, and $R_{n,2}^*$ at the 5\% significance level\label{power}}
 \renewcommand{\arraystretch}{.8}
\begin{center}
\begin{tabular}{rr|cccc|cccc}
\hline
$\sigma$&$\mu$&\multicolumn{4}{c|}{$n=25$}&\multicolumn{4}{c}{$n=75$}\\
&&$R_{n,1}$&$R_{n,2}$&$R_{n,1}^*$&$R_{n,2}^*$&$R_{n,1}$&$R_{n,2}$&$R_{n,1}^*$&$R_{n,2}^*$\\
\hline
&&\multicolumn{8}{c}{$\rho=0$}\\
1.0&1.0&28.1&18.3&8.3&6.3&57.6&41.8&11.2&8.0\\
1.0&1.5&67.0&49.7&19.2&11.3&97.5&93.0&40.2&24.8\\
0.5&1.0&46.9&85.2&12.3&70.5&88.2&99.9&21.7&99.6\\
0.5&1.5&92.2&99.2&39.2&90.6&100.0&100.0&79.7&100.0\\
\hline
&&\multicolumn{8}{c}{$\rho=0.25$}\\
1.0&1.0&7.2&6.2&10.4&7.2&6.7&6.0&16.7&10.5\\
1.0&1.5&38.8&27.0&29.6&17.5&70.9&56.9&63.9&44.8\\
0.5&1.0&22.4&77.3&16.4&78.2&43.2&99.8&32.5&99.9\\
0.5&1.5&80.9&98.5&54.0&95.5&99.7&100.0&93.5&100.0\\
\hline
&&\multicolumn{8}{c}{$\rho=0.5$}\\
1.0&1.0&1.0&1.9&15.8&9.8&0.1&1.0&32.8&20.0\\
1.0&1.5&17.7&13.1&54.7&34.6&22.4&16.6&93.7&83.2\\
0.5&1.0&8.4&71.8&24.3&91.3&7.6&99.6&53.6&100.0\\
0.5&1.5&66.0&98.1&76.4&99.5&95.7&100.0&99.5&100.0\\
\hline
&&\multicolumn{8}{c}{$\rho=-0.25$}\\
1.0&1.0&65.1&45.6&7.3&5.9&97.7&93.1&9.0&6.8\\
1.0&1.5&90.0&76.1&14.2&9.0&100.0&99.9&27.1&16.5\\
0.5&1.0&75.7&92.1&10.2&68.3&99.7&100.0&16.6&99.5\\
0.5&1.5&97.9&99.7&29.5&87.8&100.0&100.0&64.5&100.0\\
\hline
&&\multicolumn{8}{c}{$\rho=-0.5$}\\
1.0&1.0&93.8&81.0&6.7&5.7&100.0&100.0&8.1&6.4\\
1.0&1.5&99.0&94.3&11.3&7.9&100.0&100.0&19.8&12.2\\
0.5&1.0&94.9&97.8&9.0&73.9&100.0&100.0&13.3&99.8\\
0.5&1.5&99.7&100.0&23.3&90.6&100.0&100.0&50.3&100.0\\
\hline
\end{tabular}
\end{center}
\end{table}

\section{Real-Data Examples} \label{section-real-data}

\subsection{Data from karyotype analysis}
This example considers 40 unordered pairs of the lengths of the longer and shorter arms
of chromosome II of Larix decidua from 40 specimens; so $n=40$.
The data are available in Table 1 of \cite{matern1968}.
The test results from $R_{n,1}$, $R_{n,2}$, $R_{n,1}^*$, and $R_{n,2}^*$
for (\ref{hypothesis}) are as follows:
\begin{itemize}
\item
$R_{n,1}=14.91$ and $R_{n,2}=17.71$.
Calibrated by the adjusted limiting distributions,
the asymptotic $p$-values of $R_{n,1}$ and  $R_{n,2}$
are $7\times 10^{-5}$ and $2\times 10^{-4}$.

\item
$R_{n,1}^* = 1.08$ and  $R_{n,2}^* = 16.69$.
Calibrated by the adjusted limiting distributions,
the asymptotic $p$-values of $R_{n,1}^*$ and  $R_{n,2}^*$ are $0.21$ and $4\times 10^{-4}$.
 \end{itemize}
The maximum likelihood estimate  of $(\mu_1,\mu_2,\sigma_1,\sigma_2,\rho)$ is found to be
$$
(\hat\mu_1^*,\hat\mu_2^*,\hat\sigma_1^*,\hat\sigma_2^*,\hat\rho^*)
=
(62.05, 65.55, 3.50, 8.20,-0.73).
$$
Note that $\hat\rho^*=-0.73$ suggests strong negative correlation
between $X_{1i}$ and $X_{2i}$.
As revealed in the simulation studies reported in the bottom section of Table \ref{type1},
$R_{n,1}$ and $R_{n,2}$ are therefore not reliable because they are designed for $\rho = 0$.
Moreover, the fitted values $\hat\mu_1^*$ and $\hat\mu_2^*$ are very close,
but $\hat\sigma_1^*$ and $\hat\sigma_2^*$ are significantly different.
Hence, $R_{n,1}^*$ is unsuitable because it is designed for the case where $\sigma_1 = \sigma_2$. 
We recommend $R_{n,2}^*$, which
is designed to detect departures from either equal-mean or equal-variance hypotheses.


\subsection{C-band area of human chromosome data}

This example consists of normalized measurements of the C-band area on the
No.~9 chromosome pair \citep{mason1975}.
The measurements are based on three groups: the father, mother, and offspring.
These groups respectively have 40, 18, and 31 unordered pairs of normalized measurements of
the C-band area.
The data are available in Table 1 of \cite{lauder1977}.
We analyze the group of fathers as an example; the analysis of the
other groups is similar.
We constructed $R_{n,1}$, $R_{n,2}$, $R_{n,1}^*$, and $R_{n,2}^*$
and the corresponding $p$-values from the adjusted limiting distributions.
The results are as follows:
\begin{itemize}
\item
$R_{n,1}=6.51$ and $R_{n,2}=9.47$ with $n=40$.
Calibrated by the adjusted limiting distributions,
the asymptotic $p$-values of $R_{n,1}$ and  $R_{n,2}$
are $6.6\times 10^{-3}$ and $8.9\times 10^{-3}$.

\item
$R_{n,1}^* = 10.74$ and  $R_{n,2}^* = 13.48$ with $n=40$.
Calibrated by the adjusted limiting distributions,
the asymptotic $p$-values of $R_{n,1}^*$ and  $R_{n,2}^*$ are 7.5$\times10^{-4}$
and $1.9\times 10^{-3}$.
\end{itemize}
The maximum likelihood estimate  of $(\mu_1,\mu_2,\sigma_1,\sigma_2,\rho)$
is found to be
$$
(\hat\mu_1^*,\hat\mu_2^*,\hat\sigma_1^*,\hat\sigma_2^*,\hat\rho^*)
=(86.75, 68.58,10.55,8.29,0.46).
$$
Note that $\hat\rho^*=0.46$ suggests strong postive correlation
between $X_{1i}$ and $X_{2i}$.
Moreover, $\hat\mu_1^*$ and $\hat\mu_2^*$ are quite different
whereas $\hat\sigma_1^* \approx \hat\sigma_2^*$.
These suggest that $R_{n,1}^*$ is the most suitable test
while $R^*_{n,2}$ is also a possibility.
Note that $R_{n,1}^*$ is sharper than $R_{n,2}^*$ with a
smaller p-value.

\section{Technical Details\label{appendix}}
\subsection{Reparameterization and preparation lemmas}

Recall that $(Y_{1i}, Y_{2i})$ is the unordered pair of $(X_{1i}, X_{2i})$
and the latter has a bivariate normal distribution with parameter vector
$\btheta = (\mu_1, \mu_2, \sigma_1, \sigma_2, \rho)^\tau$.
The log-likelihood function based on $\{(Y_{1i}, Y_{2i})\}_{i=1}^n$ is
\begin{eqnarray*}
\ell_n^*(\btheta)
&=&
\sum_{i=1}^{n}\log\{\phi(Y_{1i}, Y_{2i}; \btheta)+\phi(Y_{2i}, Y_{1i};\btheta)\}\\
&=&
\sum_{i=1}^{n} \log\{\phi (X_{1i},  X_{2i}; \btheta)+\phi (X_{2i}, X_{1i}; \btheta)\}.
\end{eqnarray*}

Let $Z_{1i}=(X_{1i}+X_{2i})/2$ and $Z_{2i}=(X_{1i}-X_{2i})/2$.
We introduce notation for the following quantities:
\begin{align*}
&
\E (Z_{1i}) = (\mu_1+\mu_2)/2 = \mu, \\
&
 \E(Z_{2i}) = (\mu_1-\mu_2)/2 = \Delta,\\
&\var(Z_{1i})
 = (1/4)(\sigma_1^2+\sigma_2^2+2\rho\sigma_1\sigma_2)= \sigma^2_+,\\
&\var(Z_{2i})
 = (1/4)(\sigma_1^2+\sigma_2^2- 2\rho\sigma_1\sigma_2) = \sigma^2_-,\\
&\cov(Z_{1i}, Z_{2i})
 = (1/4) (\sigma_1^2-\sigma_2^2) = \xi \sigma_+\sigma_-.
\end{align*}
Further, let
\(
\beta_0 = \Delta - \mu  ({\sigma_-}/{\sigma_+}) \xi,~~
\beta_1 = ({\sigma_-}/{\sigma_+}) \xi,~~
\eta^2 = (1-\xi^2)\sigma_-^2,
\)
and
\begin{eqnarray*}
\ell_{n,1}^*(\mu, \sigma_+)
&=&
\sum_{i=1}^{n} \log\{\phi (Z_{1i}; \mu, \sigma_+)\},
\\
\ell_{n,2}^*(\beta_0, \beta_1, \eta)
&=&
\sum_{i=1}^{n}
\log
\{
0.5\phi (Z_{2i}; \beta_0+\beta_1Z_{1i}, \eta)
+
0.5\phi (- Z_{2i}; \beta_0+ \beta_1Z_{1i}, \eta)
\}.
\end{eqnarray*}
Note that we use $\phi(x; \mu,\sigma)$ to denote the density function of $N(\mu, \sigma^2)$,
matching $\phi(x_1, x_2; \btheta)$ for the bivariate normal distribution.

With these, we obtain the following decomposition of the likelihood function:
\[
\ell_n^*(\btheta)
=
\ell_{n,1}^*(\mu, \sigma_+)
+
\ell_{n, 2}^*(\beta_0, \beta_1,\eta).
\]
We use a generic $\btheta$ for the parameters,
which may be interpreted as $\btheta = (\mu, \sigma_+, \beta_0, \beta_1, \eta)^\tau$
when necessary.

Under $H_0$ in Theorem \ref{theorem-1} which includes the
assumption that $\rho = 0$, suppose the true parameter values of the
data-generating distribution are $\mu_1 = \mu_2 = \mu_*$,
$\sigma^2_1 = \sigma^2_2 = \sigma^2_*$.
We may then, in our proofs, work with the transformed data
\[
X_1^* = \sqrt{2}(X_1 - \mu_*)/\sigma_*, ~~
X_2^* = \sqrt{2}(X_2 - \mu_*)/\sigma_*.
\]
After the transformation, the algebraic form of the likelihood does not change
but the true parameter values of the data-generating distribution
become $\mu_1 = \mu_2 = 0$ and $\sigma_1^2 = \sigma_2^2 = 2$.
{Without loss of generality,  based on the above invariance property,
we may assume that the true parameters $\mu_1 = \mu_2 = 0$ and
$\sigma_1^2 = \sigma_2^2 = 2$ under $H_0$.}

Under $H_0$ in Theorem \ref{theorem-2}, without loss of
generality, the same assumption is applicable to $\mu$ and $\sigma$.
We now reveal that by the same invariance principle we may also assume $\rho =0$ as
long as the true value $\rho \neq \pm 1$.
When $\rho_* \neq \pm 1$, we simply let
\[
(X_1^{**}, X_2^{**}) = \{X^*_1 , ~ (X^*_2 - \rho_* X^*_1)/\sqrt{1 - \rho_*^2} \}.
\]
The distribution-generated data $\{X_1^{**}, X_2^{**}\}$ now  has
the true parameter values $\mu_1 = \mu_2 = 0$, $\sigma_1^2 = \sigma_2^2 = 2$,
and $\rho = 0$ under $H_0$.

With the above standardization operation, for both Theorems \ref{theorem-1} and \ref{theorem-2},
we study the asymptotic null properties under the assumption that
$Z_{1i}$ and $Z_{2i}$ are independent normal random variables
with the standard parameter values:
\[
(\mu, \sigma_+, \beta_0, \beta_1, \eta) = (0, 1, 0, 0, 1).
\]
We first establish three preparatory lemmas.

\begin{lemma}
\label{lemma0}
As $n \to \infty$, we have, almost surely,
\[
\sup_{\beta_0, \beta_1} \sum_{i=1}^n \ind (|Z_{2i}-\beta_0-\beta_1Z_{1i}| \leq 1/4 )\leq (1/4) n,
\]
where $\ind(\cdot)$ is the indicator function.
\end{lemma}

\proof
Note that
\[
n^{-1} \sum_{i=1}^n \ind (|Z_{2i}-\beta_0-\beta_1Z_{1i}| \leq 1/4 )
\]
is the empirical measure of the two-dimensional stripe formed by
the inequality
\[
|Z_{2}-\beta_0-\beta_1Z_{1}| \leq 1/4.
\]
This class of stripes can divide $n$ points in two-dimensional space
into at most a polynomial number of different subsets.
By \cite{Pollard1990}, this property
 implies the uniform strong law of large numbers:
\begin{equation}
\label{lemma0.eq3}
\sup_{\beta_0,\beta1}
\left|
n^{-1}
\sum_{i=1}^n
\ind (|Z_{2i}-\beta_0-\beta_1Z_{1i}|\leq 1/4)
-P\big (|Z_2-\beta_0-\beta_1Z_1|\leq 1/4 \big )
\right|
\to 0
\end{equation}
almost surely.

The distribution of $Z_{2}-\beta_0-\beta_1Z_{1}$ is normal with
{variance at least 1.}
Based on this, we have
$P( |Z_{2}-\beta_0-\beta_1Z_{1}| \leq1/4) \leq {0.2}$
for any $\beta_0, \beta_1$.
Hence, almost surely,
\[
\sum_{i=1}^n
\ind (|Z_{2i}-\beta_0-\beta_1Z_{1i}|\leq 1/4)
\leq {0.2} n + o(n)
\leq 0.25n.
\]
This completes the proof.
\qed


\begin{lemma}
\label{lemma1}
Suppose an estimator $\bar \btheta$ satisfies
\begin{eqnarray}
\ell_n(\bar \btheta) - \ell_n(\btheta_0)
&=&
\{\ell_{n,1}^*(\bar \mu,\bar \sigma_+)
+ \ell_{n,2}^*(\bar\beta_0,\bar\beta_1, \bar\eta)\}
- \{\ell_{n,1}^*(0,1)\} + \ell_{n,2}^*(0,0,1)\}
\nonumber \\
&=&
\{\ell_{n,1}^*(\bar \mu,\bar \sigma_+)  - \ell_{n,1}^*(0,1)\}
+
\{\ell_{n,2}^*(\bar\beta_0,\bar\beta_1, \bar\eta) - \ell_{n,2}^*(0,0,1)\}
\nonumber \\
& \geq&
C > -\infty
\label{condition.inequality}
\end{eqnarray}
for some constant $C$.
Then under the null model, $\bar\btheta= \btheta_0 + o_p(1) = (0,1,0,0,1)^\tau+o_p(1)$.
\end{lemma}

\proof
Note that we have decomposed $\ell_n(\bar \btheta) - \ell_n(\btheta_0)$ into
a sum of two terms. For the first term, according to the classical result about the
LRT under regular models, it is clear that
\begin{equation}
\label{chen1}
\sup_{\mu, \sigma_+}  \{\ell_{n,1}^*(\mu, \sigma_+)  - \ell_{n,1}^*(0,1)\}
= O_p(1).
\end{equation}

When in the second term the variance parameter $\eta > M_0 = \exp(4)$, we have
\[
\sum_{i=1}^{n}
\log
\{
0.5\phi (Z_{2i}; \beta_0+\beta_1Z_{1i}, \eta)
+
0.5\phi (- Z_{2i}; \beta_0+ \beta_1Z_{1i}, \eta)
\}
\leq
- n \log M_0 = - 4 n.
\]
By the law of large numbers, we have
\[
n^{-1} \ell_{n, 2}^*(0, 0, 1) \geq -(1/2)\log (2 \pi) - \E(Z_2^2) \geq -2.
\]
almost surely. This implies that
\[
\ell_{n, 2}^*(\beta_0, \beta_1, \eta) - \ell_{n, 2}^*(0, 0, 1)
\leq  - 2n
\]
and subsequently, uniformly for $\eta$ in this range,
\[
\ell_{n, 2}^*(\beta_0, \beta_1,\eta) -  \ell_{n, 2}^*(0, 0, 1) \to - \infty.
\]
Together with \eqref{chen1}, we have, whenever $\eta > M_0 = \exp(4)$,
\[
\ell_n(\btheta) - \ell_n(\btheta_0) \to - \infty
\]
in probability.
Since the lemma condition clearly states that $\bar \eta$ does not
have the above property, it cannot be in this range. That is,
we conclude that $\bar \eta \leq M_0$.

Suppose $\eta < \epsilon_0$ and $\epsilon_0$ is a very small positive value.
In this case, for all $i$, we have
\[
\log
\{
0.5\phi (Z_{2i}; \beta_0+\beta_1Z_{1i}, \eta)
+
0.5\phi (- Z_{2i}; \beta_0+ \beta_1Z_{1i}, \eta)
\}
\leq
-  \log (\eta).
\]
For $i$ such that
\begin{equation}
\label{chen2}
\min \{ |Z_{2i}+ \beta_0+\beta_1Z_{1i}|, |Z_{2i}- \beta_0- \beta_1Z_{1i}|\} > 1/4,
\end{equation}
we have
\[
\log
\{
0.5\phi (Z_{2i}; \beta_0+\beta_1Z_{1i}, \eta)
+
0.5\phi (- Z_{2i}; \beta_0+ \beta_1Z_{1i}, \eta)
\}
\leq
- \log (\eta) - (1/32)/\eta^2.
\]
By Lemma \ref{lemma0}, uniformly in $\beta_0$ and $\beta_1$ and almost surely,
at least $(1/2)n$ of the $i$'s satisfy \eqref{chen2}.
Therefore,
\[
\ell_{n, 2}^*(\beta_0, \beta_1,\eta) - \ell_{n, 2}^*(0, 0, 1)
\leq
- \{ \log (\eta)  + (1/64)/\eta^2\} n \to - \infty
\]
as $n \to \infty$ and $\eta \to 0$.
Namely, for all $\eta < \epsilon_0$ sufficiently small, we also have
\[
\ell_{n, 2}^*(\beta_0, \beta_1,\eta) -  \ell_{n, 2}^*(0, 0, 1) \to - \infty.
\]
In conclusion, the $\bar \eta$ value satisfying the lemma condition
must almost surely fall within the interval $[\epsilon_0, M_0]$ for some
sufficiently small $\epsilon_0 > 0$ and sufficiently large $M_0 < \infty$.

Within the parameter space $[\epsilon_0, M_0] \times \mathbb{R}^2$,
the density function
\[
0.5\phi (Z_{2i}; \beta_0+\beta_1Z_{1i}, \eta)
+
0.5\phi (- Z_{2i}; \beta_0+ \beta_1Z_{1i}, \eta)
\]
satisfies the conditions for the consistency of the MLE specified in
\cite{wald1949}.
For instance, it is a continuous density function with its limit being 0
whenever $\beta_0$ or $\beta_1$ goes to infinity.
For a sufficiently small $\epsilon > 0$, let
\[
B_{\epsilon}
=\{(\beta_0, \beta_1, \eta):  \beta_0^2+\beta_1^2+(\eta^2-1)^2 \leq \epsilon^2\}
\]
be a ball centered at the true value.
The side conclusion as stated in \cite{wald1949} is
\begin{eqnarray}
\label{ln2.part2.upper2}
\sup_{(\beta_0,\beta_1,\eta) \not \in B_{\epsilon}}
\ell_{n,2}^*(\beta_0,\beta_1,\eta)
-
\ell_{n,2}^* (0, 0, 1) \leq - \delta n \to - \infty
\end{eqnarray}
for some $\delta > 0$.
Again, by the lemma condition on $\bar{\btheta}$, we must
have $\bar\beta_0, \bar \beta_1, \bar \eta$ within $\epsilon$
of the true parameter value for any $\epsilon > 0$ as $n \to \infty$.
This proves part of the lemma.

It is now apparent that we also have
\[
\sup_{\beta_0, \beta_1,\eta} \{ \ell_{n,2}^*(\beta_0,\beta_1,\eta)
-
\ell_{n,2}^* (0, 0, 1)\}
= O_p(1).
\]
By the same argument based on the assumed property of $\bar \btheta$,
we must have
\[
\ell_{n,1}^*(\bar \mu, \bar \sigma_+)  - \ell_{n,1}^*(0,1)
= O_p(1) = o_p(n).
\]
This is sufficient for the proof of the consistency of $(\bar \mu, \bar \sigma_+)$.
Combined with the proof of the other parts,
this completes the proof of the lemma.
\qed

Next, we strengthen the results of Lemma \ref{lemma1}.
We first define some notation for the next lemma.
Let
\begin{align*}
A_i &=
(Z_{1i}, {(Z_{1i}^2-1)}/{2})^\tau,\\
B_i & =
((Z_{2i}^2-1)/2, (Z_{1i}^2-1)(Z_{2i}^2-1)/2,
Z_{1i}(Z_{2i}^2-1)/2, - (Z_{2i}^4-6Z_{2i}^2+3)/12)^\tau.
\end{align*}
It can be seen that $\E(A_i)= 0$, $\E(B_i) = 0$, $A_i$ and $B_i$ are
uncorrelated, and
\[
\Sigma_A = \var (A_i) =  \mbox{diag}(1,1/2); ~~~
\Sigma_B = \var (B_i) = \mbox{diag}(1/2, 1, 2, 1/6).
\]
Further, we introduce two parameter vectors of lengths 2 and 4:
\[
\ss_1
=(\mu, ~\sigma_+^2-1)^\tau; ~~
\ss_2  =
(\beta_0^2+\beta_1^2 + (\eta^2-1),~ \beta_1^2,~\beta_0\beta_1,~\beta_0^4)^\tau.
\]
{In the following, we use
$|{\bf x}|$ and $\|{\bf x}\|$
to denote the $L_1$ and $L_2$ norms of the vector ${\bf x}$, respectively. }

\begin{lemma}
\label{lemma2}
Under the conditions of Lemma \ref{lemma1}
and the null hypothesis, we have
\begin{align*}
\hspace{-1em}
(a)~~ &
\ell_{n,1}^*(\bar\mu, \bar\sigma_+) - \ell_{n,1}^*(0,1)
=
 \bar \ss_1^\tau \sum_{i=1}^n A_i - (n/2) \{ \bar \ss_1^\tau \Sigma_A \bar \ss_1 \} \{1+o_p(1)\}+o_p(1);
\vspace{0.5ex} \\
\hspace{-1em}
(b)~~ &
\ell_{n,2}^*(\bar\beta_0,\bar\beta_1,\bar\eta) - \ell_{n,2}^*(0,0,1)
\leq
\bar \ss_2^\tau \sum_{i=1}^n B_i - (n/2) \{ \bar \ss_2^\tau \Sigma_B \bar \ss_2 \} \{1+o_p(1)\}+o_p(1);
\\
\hspace{-1em}
(c)~~ &
\bar\mu, ~\bar\sigma_+^2-1, ~ \bar\beta_0^4, ~\bar\beta^2_1
\mbox{~and}
 ~(\bar\eta^2-1)^2 \mbox{~are~} O_p(n^{-1/2}).
\end{align*}
\end{lemma}
\proof

We first prove (a). By Lemma \ref{lemma1},
we have
$
(\bar\mu, \bar\sigma_+)=(0,1)+o_p(1)
$.
We obtain (a) by expanding
$\ell_{n,1}^*(\bar\mu, \bar\sigma_+)$ at $(\bar\mu, \bar\sigma_+)=(0,1)$
to the second order and then assessing the asymptotic orders via
the weak law of large numbers.

To prove (b), we first denote
 $$
 \delta_i (\beta_0,\beta_1,\eta)
 =
 \{\phi (Z_{2i};\beta_0+\beta_1Z_{1i},\eta) + \phi (-Z_{2i}; \beta_0+\beta_1Z_{1i},\eta)\}/
 \{2\phi(Z_{2i}; 0, 1)\} - 1
 $$
 and then write
 $$
\ell_{n,2}^*(\beta_0,\beta_1,\eta) - \ell_{n,2}^*(0,0,1)
 =\sum_{i=1}^n \log\{ 1+\delta_i (\beta_0,\beta_1,\eta)\}.
 $$
Applying the inequality $\log(1+x)\leq x-x^2/2+x^3/3$, we have
\be
\label{Rn.star.21}
\ell_{n,2}^*(\beta_0,\beta_1,\eta) - \ell_{n,2}^*(0,0,1)
 \leq
\sum_{i=1}^n \delta_i(\beta_0,\beta_1,\eta)
-
(1/2)\sum_{i=1}^n \delta_i^2(\beta_0,\beta_1,\eta)
+
(1/3) \sum_{i=1}^n \delta_i^3(\beta_0,\beta_1,\eta).
\ee
Next, we delineate $\delta_i(\beta_0, \beta_1, \eta)$ given
$
(\bar\beta_0, \bar\beta_1,\bar\eta) =(0, 0, 1)+o_p(1)
$
as proved in Lemma \ref{lemma1}.
We perform two main steps.
In the {\it first step}, we obtain the fourth-order Taylor expansion of
$\delta_i(\beta_0, \beta_1, \eta)$;
in the {\it second step}, we assess the asymptotic orders of
the terms in the expansion and put them into appropriate order expressions.

We start with the {\it first step}.
Let the partial derivatives be
$$
\delta^{(s,t,k)}_i(\beta_0,\beta_1,\eta)
=\frac{\partial ^{s+t+k}\delta_i(\beta_0,\beta_1,\eta) }
        {\partial \beta_0^s\partial\beta_1^t\partial (\eta^2)^k}.
$$
Expanding both
$\phi(\pm Z_{2i}; \beta_0+\beta_1Z_{1i},\eta)$
to the fourth order at $(\beta_0,\beta_1,\eta) = (0,0,1)$,
we get
\begin{equation}
\label{deltai1}
\delta_i(\beta_0, \beta_1,\eta)
=
\sum_{s+t+k=1}^{4}\frac{\beta_0^s\beta_1^t(\eta^2-1)^k}
{s!t!k!}\delta^{(s,t,k)}_i (0,0,1)+\epsilon_{in}^{(1)},
\end{equation}
where the summation is over all non-negative integer combinations of
$s, t, k$ summing to $4$ and
$\epsilon_{in}^{(1)}$ is the remainder term in the Taylor expansion.
Let $\epsilon_n^{(1)}=\sum_{i=1}^n\epsilon_{in}^{(1)}$, then
\begin{equation*}
\label{epsiloni1}
\epsilon_n^{(1)}
=O_p(n^{1/2})\sum_{s+t+k=5} \beta_0^s\beta_1^t(\eta^2-1)^k
=o_p(n^{1/2}){|\ss_2|}.
\end{equation*}

In the {\it second step}, we first show that every term in the summation part of
(\ref{deltai1})
satisfying $s+2t+2k\geq 5$ is of order {$o_p(n^{1/2}){| \ss_2 |}$}.
For instance, when $s=t=k= 1$, we have
\[
|\beta_0\beta_1(\eta^2-1)|
\leq
|\beta_0| \{\beta_1^2+ (\eta^2-1)^2\}
= o_p(|\ss_2|),
\]
helped by the fact that we are investigating the region of $\beta_0 = o_p(1)$.
For notational simplicity, let $\delta^{(s,t,k)}_i= \delta^{(s,t,k)}_i( 0, 0, 1)$.
It is easy to check that $\delta^{(s,t,k)}_i$ has zero mean and finite variance, so
\[
\sum_{i=1}^n \delta^{(s,t,k)}_i= O_p(n^{1/2}).
\]
Therefore, we have
$$
\sum_{i=1}^n
\frac{\beta_0^s\beta_1^t(\eta^2-1)^k} {s!t!k!}\delta^{(s,t,k)}_i
=
o_p(n^{1/2} | \ss_2|).
$$
The proofs for the other $s+2t+2k\geq 5$ terms are similar.
Hence, we may write
\begin{equation}
\label{deltai2}
\delta_i(\beta_0,\beta_1,\eta)
=
\sum_{s+2t+2k=1}^{4}
\frac{\beta_0^s\beta_1^t(\eta^2-1)^k} {s!t!k!}\delta^{(s,t,k)}_i +\epsilon_{in}^{(2)}
\end{equation}
and still have
\begin{equation}
\label{epsiloni2}
\sum_{i=1}^n\epsilon_{in}^{(2)}
=
o_p(n^{1/2} |\ss_2|).
\end{equation}
By straightforward algebra, we find
\be
\sum_{s+2t+2k=1}^{4}
\frac{\beta_0^s\beta_1^t(\eta^2-1)^k} {s!t!k!}\delta^{(s,t,k)}_i
=
\ss_2^\tau B_i -1.5\{\beta_0^2 + (\eta^2-1)\}^2B_{i}[4]
\ee
where the unwanted term $B_{i}[4]$ is the fourth element of vector $B_i$.
Its coefficient is easily verified to be $\{\beta_0^2 + (\eta^2-1)\}^2 = o_p(|\ss_2|)$.
This allows us to obtain a neater expression by absorbing it into the higher-order term,
concluding that
\begin{equation}
\label{deltai3}
\delta_i(\beta_0,\beta_1,\eta)
=
\ss_2^\tau B_i+\epsilon_{in}^{(3)}
\end{equation}
such that
\begin{equation}
\label{epsiloni3}
\sum_{i=1}^n\epsilon_{in}^{(3)}
= o_p(n^{1/2} |\ss_2|) = o_p(1)+o_p(n\|\ss_2\|^2).
\end{equation}
In short, we have shown that
\begin{eqnarray}
 \sum_{i=1}^n\delta_i(\beta_0,\beta_1,\eta)
 =
\ss_2^\tau  \sum_{i=1}^n B_i+o_p(1) + o_p(n\|\ss_2\|^2).
 \label{Rn.star.21.linear}
\end{eqnarray}

The above algebraic manipulations are typical of the techniques
employed in \cite{chen2001} and \cite{chenjrss2001}.
The same techniques, which are tedious but not sophisticated,
give
\begin{eqnarray*}
\sum_{i=1}^n \delta^2_i(\beta_0, \beta_1, \eta)
&=&
\ss_2^\tau  \big \{ \sum_{i=1}^n B_iB_i^\tau \big \} \ss_2 +o_p(1) + o_p(n\|\ss_2\|^2),
\\
\sum_{i=1}^n \delta^3_i(\beta_0,\beta_1,\eta)
&=&
\sum_{i=1}^n |\ss_2^\tau B_i| ^3 + o_p(1) + o_p(n\|\ss_2\|^2).
\end{eqnarray*}
Together with the weak law of large numbers these lead to
\begin{eqnarray}
\sum_{i=1}^n \delta^2_i(\beta_0, \beta_1, \eta)
&=&
n \ss_2^\tau \Sigma_B \ss_2 +o_p(1)+ o_p(n\|\ss_2\|^2),
 \label{Rn.star.21.quad}
\\
\sum_{i=1}^n \delta^3_i(\beta_0,\beta_1,\eta)
 &=&
 o_p(1) + o_p(n\|\ss_2\|^2).
 \label{Rn.star.21.cubic}
\end{eqnarray}
Combining (\ref{Rn.star.21.linear})--(\ref{Rn.star.21.cubic}) with (\ref{Rn.star.21}),
we have
\[
\ell_{n,2}^*(\beta_0, \beta_1, \eta) - \ell_{n,2}^*(0,0,1)
 \leq
\ss_2^\tau  \sum_{i=1}^n B_i
 - (n/2) \ss_2^\tau \Sigma_B \ss_2 \{1+o_p(1)\} + o_p(1).
 \label{Rn.star.21.upper}
\]
Recall that $(\bar\beta_0,\bar\beta_1,\bar\eta) = (0, 0, 1) + o_p(1)$,
so the above upper bound is applicable to
$\ell_{n,2}^*(\bar \beta_0, \bar \beta_1, \bar\eta) - \ell_{n,2}^*(0,0,1)$.
This completes the proof of (b).

Finally, we come to (c).
Combining (a) and (b) and the conditions in Lemma \ref{lemma1},
 we have
\begin{eqnarray}
C
&\leq&
\{\ell_{n,1}^*(\bar \mu,\bar \sigma_+)- \ell_{n,1}^*(0,1)\}
+\{\ell_{n,2}^*(\bar\beta_0, \bar\beta_1, \bar\eta)- \ell_{n,2}^*(0,0,1)\}
\nonumber\\
&\leq&
\sum_{i=1}^n\{ \bar \ss_1^\tau A_i + \bar \ss_2^\tau B_i  \}
- (n/2) \{\bar\ss_1^\tau \Sigma_A \bar\ss_1 + \bar\ss_2^\tau \Sigma_B \bar\ss_2 \} \{1+o_p(1)\}+o_p(1),
 \label{lik.upper1}
\end{eqnarray}
which is possible only if both
$
\bar \ss_1 =O_p(n^{-1/2})$ and $\bar \ss_2 = O_p(n^{-1/2})
$.
This leads to the order assessments in (c)
and completes the proof of the entire lemma.
\qed

\subsection{Proof of Theorem \ref{theorem-1}}
The difference between Theorems \ref{theorem-1} and \ref{theorem-2}
is that in the former we consider $\rho_0 = 0$ to be known when formulating the test statistic.
This makes it helpful to reorganize the entries of $A_i$ and $B_i$
and the corresponding entries of $\ss_1$ and $\ss_2$.

When $\rho_0 = 0$ is known, we have $\sigma_+ = \sigma_-$.
Let
\[
\tt = (\mu, \beta_0^2/2 + \sigma_+^2 -1, \beta_0^2, \beta_1^2, \beta_0 \beta_1)^\tau.
\]
Every entry of $\ss_1$ and $\ss_2$ is a linear combination
of the entries of $\tt$, possibly with an $O_p(\| \tt \|^2)$ difference when these parameter
values approach their default null values.
We enumerate these entries as follows. The first entry of $\ss_1$ is  $\ss_1[1] = \tt[1]$,
and the second is $\ss_1[2] = \tt[2]-\tt[3]/2$.
For the entries of $\ss_2$, we have
\[
\ss_2[1] = \beta_0^2 + \beta_1^2 + (\eta^2 - 1)=
   \tt[2] +\tt[3]/2- \beta_1^2(\sigma^2_+-1) = \tt[2] +\tt[3]/2+ O_p( \|\tt\|^2).
\]
For the others, $\ss_2[2] = \tt[4]$, $\ss_2[3] = \tt[5]$,
and $\ss_2[4] = (\tt[3] )^2= O_p(\| \tt \|^2)$.

Because every entry of $\ss_1$ and $\ss_2$ is virtually a linear combination
of the entries of $\tt$, we can reorganize the entries of $A_i$ and $B_i$ into
a vector $D_i$ such that
\[
\ss_1^\tau A_i + \ss_2^\tau B_i  = \{\tt^\tau + O_p(\| \tt\|^2)\} D_i.
\]
Naturally, we have $\E( D_i) = 0$ and some algebra shows that
$\var(D_i) = \Sigma_D = \mbox{diag}(1, 1, 1/4, 1, 2)$.
The following result is immediate.

\begin{lemma}
\label{lemma3}
Assume the conditions of Lemma \ref{lemma2} and let $\bar\rho=0$.
If, under the null model,
$$
\ell_n(\bar \btheta) - \ell_n(\btheta_0) \geq C>-\infty,
$$
we then have
\begin{align*}
(a) ~~&
\ell_n(\bar \btheta) - \ell_n(\btheta_0)
\leq
\bar \tt^\tau \sum_{i=1}^n D_i
- (n/2) \bar \tt^\tau \Sigma_D \bar\tt  \{1+o_p(1)\}+o_p(1); \\
(b) ~~&
 \bar \mu, ~\bar\sigma_+^2-1, ~\bar\beta_0^2,
 \mbox{~and~} \bar\beta_1^2 \mbox{ are } O_p(n^{-1/2}).
\end{align*}
\end{lemma}

We are now ready for Theorem \ref{theorem-1}.
The order conclusions of the MLEs in both
Theorem \ref{theorem-1}(a) and \ref{theorem-1}(b) have been
established in Lemma \ref{lemma3}.
We now derive the limiting distributions.

We rewrite  $R_{n, 1}$ defined in \eqref{Rn1} as
\[
R_{n, 1} = 2 \{\ell_n(\tilde \btheta ) - \ell_n(\btheta_0 )\}
- 2  \{\ell_n(\check \btheta ) - \ell_n(\btheta_0 ) \}
\]
with $\check \btheta$ being the maximum point of the reduced model
where $(\mu_1, \sigma_1) = (\mu_2,\sigma_2)$.
Since the reduced model is regular, by standard techniques such as
those in \citet{serfling2000}:
\begin{equation}
\label{lik.dif.01}
2\{\ell_n( \check \btheta ) - \ell_n(\btheta_0)\}
=
n^{-1}\Big \{
\big (\sum_{i=1}^n D_i[1] \big )^2  + \big (\sum_{i=1}^n D_i[2] \big )^2
\Big \}
+o_p(1)
\end{equation}
where $D_i[1], D_i[2]$ denote the first two entries of vector $D_i$.

Next, note that $\tilde \btheta$ is the maximum point of the reduced model
where $\sigma_1=\sigma_2=\sigma$.
This makes $\beta_1=\xi=0$ and subsequently for $\tt$ under the reduced model,
\[
\tt = (\mu, \beta_0^2/2 + (\sigma_+^2 -1), \beta_0^2, 0, 0)^\tau.
\]
Nevertheless, Lemma \ref{lemma3} is applicable to the above form of
$\tt$ as long as it is close to its counterpart in the null model. Hence,
\begin{eqnarray}
2\{ \ell_n(\btheta) - \ell_n(\btheta_0)\}
&\leq&
2 \sum_{i=1}^n \tt^\tau D_i
- n \tt^\tau \Sigma_D \tt
+o_p(1)
\nonumber\\
&\leq&
\sup
\left \{
2\sum_{i=1}^n \tt^\tau D_i
- n  \tt^\tau\Sigma_D  \tt : \tt[3]\geq0, \tt[4]=0, \tt[5]=0
\right \}
+o_p(1)
\nonumber\\
&\leq&
n^{-1}
\left [
(\sum_{i=1}^n D_i[1])^2
+
(\sum_{i=1}^n D_i[2])^2
+
 4\{\big (\sum_{i=1}^n D_i[3] \big )^{+}\}^2\right] +o_p(1).
\label{lik.dif.11.upper}
\end{eqnarray}
Note the range of the supremum conforms to the form of $\tt$ in the
reduced model and the fact that $\tt[3] =\beta_0^2\geq 0$.
The specific coefficient values are due to the value of $\Sigma_D$.

The upper bound in (\ref{lik.dif.11.upper})  is attained if we put
\[
\tt = n^{-1}
\Big (
\sum_{i=1}^n  D_i[1],~~\sum_{i=1}^n  D_i[2],~~4 \big \{ \sum_{i=1}^n D_i[3]\big \}^{+},
0, 0
\Big )^\tau.
\]
With some straightforward algebra, the corresponding $\btheta$ values of $\tt$ exist and satisfy
\begin{eqnarray}
\mu_1=O_p(n^{-1/4}),~~
\mu_2 =O_p(n^{-1/4}), ~~
\sigma^2-1=O_p(n^{-1/2}).
\end{eqnarray}
Applying the Taylor expansion,  with $\btheta$ being the above $\btheta$,
we get
\be
2\{\ell_n(\btheta) - \ell_n(\btheta_0)\}
=
n^{-1}
\left [
(\sum_{i=1}^n D_i[1])^2
+
(\sum_{i=1}^n D_i[2])^2
+
 4\{\big (\sum_{i=1}^n D_i[3] \big )^{+}\}^2\right]  +o_p(1).
\label{lik.dif.11.lower0}
\ee
Since $\tilde \btheta$ is the maximum point of $\ell_n(\btheta)$,
$
2\{\ell_n(\tilde\btheta) - \ell_n(\btheta_0)\}
$
is not smaller than the value in  (\ref{lik.dif.11.upper}).
The sandwich technique of \cite{chen2001} and \cite{chenjrss2001} or the squeeze theorem
can be applied to obtain
\be
2\{\ell_n(\tilde\btheta) - \ell_n(\btheta_0)\}
=
n^{-1}
\left [
(\sum_{i=1}^n D_i[1])^2
+
(\sum_{i=1}^n D_i[2])^2
+
 4\{\big (\sum_{i=1}^n D_i[3] \big )^{+}\}^2\right] +o_p(1).
\label{lik.dif.11}
\ee

Combining (\ref{lik.dif.01}) and (\ref{lik.dif.11}) gives
$$
R_{n,1}=  4 n^{-1} \{\big (\sum_{i=1}^n D_i[3] \big )^{+}\}^2 +o_p(1),
$$
which has the limiting distribution $0.5\chi^2_0+0.5\chi^2_1$.
This completes the proof of part (a).

We now prove conclusion (b). In this case,
the range of $\tt$ has only an intrinsic restriction as seen in the expression
\[
\tt = (\mu, \beta_0^2/2 + (\sigma_+^2 -1), \beta_0^2, \beta_1^2, \beta_0 \beta_1)^\tau.
\]
Let
$\tt_1= (\mu, \beta_0^2/2 + (\sigma_+^2 -1))^\tau$
and $\tt_2 = (\beta_0^2, \beta_1^2, \beta_0 \beta_1)^\tau$.
It can be seen that $\tt_2$ lies on a two-dimensional manifold.
Nonetheless, the upper bound developed in Lemma \ref{lemma3} remains valid.
We partition $D_i$ into $D_{i1}$ and $D_{i2}$ with covariance matrices
$\Sigma_{D1}$ and $\Sigma_{D2}$.
With these preparations, we have
\begin{eqnarray}
2\{\ell_n(\hat \btheta) - \ell_n(\btheta_0)\}
&\leq&
2 \hat\tt_1^\tau \sum_{i=1}^n D_{i1} + 2\hat\tt_2^\tau  \sum_{i=1}^n D_{i2}
-
n \big (\hat \tt_1^\tau \Sigma_{D1} \hat\tt_1 + \hat \tt_2^\tau \Sigma_{D2}\hat \tt_2\big ) +o_p(1)
\nonumber\\
&\leq&
n^{-1} (\sum_{i=1}^n D_{i1})^\tau (\sum_{i=1}^n D_{i1})
+
\sup_{\tt_2} \big \{ 2 \tt_2^\tau  \sum_{i=1}^n D_{i2} - n\tt_2^\tau \Sigma_{D2} \tt_2\big \}
+o_p(1).
\label{lik.dif.12.upper}
\end{eqnarray}
The supremum is taken over $\tt_2$ with the intrinsic restriction respected.
Similarly to (\ref{lik.dif.11}), the upper bound (\ref{lik.dif.12.upper}) is attained
at some feasible parameter value.
Hence,
 \be
2\{\ell_n(\hat \btheta) - \ell_n(\btheta_0)\}
=
n^{-1} (\sum_{i=1}^n D_{i1})^\tau (\sum_{i=1}^n D_{i1})
+
\sup_{\tt_2} \big \{ 2 \tt_2^\tau  \sum_{i=1}^n D_{i2} - n\tt_2^\tau \Sigma_{D2} \tt_2\big \}
+o_p(1).
\label{lik.dif.12}
\ee
Combining (\ref{lik.dif.01}) and (\ref{lik.dif.12}), we get
\[
R_{n,2}
=
\sup_{\tt_2} \big \{ 2 \tt_2^\tau  \sum_{i=1}^n D_{i2} - n\tt_2^\tau \Sigma_{D2} \tt_2\big \}
+ o_p(1).
\]
The intrinsic restriction due to the specific form of $\tt_2 = (\beta_0^2, \beta_1^2, \beta_0 \beta_1)^\tau$
leads to the nonstandard form of the limiting distribution
in the theorem.

\subsection{Proof of Theorem \ref{theorem-2}}
The test problem in Theorem \ref{theorem-2} is different from that of
Theorem \ref{theorem-1} because we do not assume knowledge of the $\rho_0$ value.
The parameter vector is now $\btheta = (\mu_1,\mu_2,\sigma_1,\sigma_2,\rho)^\tau$
including the correlation coefficient $\rho$.
Because of the invariance argument, we need consider only the case where
$\btheta_0  =  (0, 0, 1, 1, 0)^\tau$ under the null hypothesis
for the asymptotic properties in this theorem.

With the introduction of $\rho$, it helps to redefine $\ss_1$, $\ss_2$, and so on
as follows:
\[
\ss_1
=(\mu, ~\sigma_+^2-1, \beta_0^2+\beta_1^2 + (\eta^2-1))^\tau; ~~
\ss_2  =
(\beta_1^2,~\beta_0\beta_1,~\beta_0^4)^\tau
\]
and the corresponding $A_i$, $B_i$ as
\begin{align*}
A_i &=
(Z_{1i}, {(Z_{1i}^2-1)}/{2}, (Z_{2i}^2-1)/2)^\tau,\\
B_i & =
((Z_{1i}^2-1)(Z_{2i}^2-1)/2,  Z_{1i}(Z_{2i}^2-1)/2, - (Z_{2i}^4-6Z_{2i}^2+3)/12)^\tau.
\end{align*}
These are almost the quantities with the same names defined above Lemma \ref{lemma2}.
The difference is that the first entry of $\ss_2$ is now the third entry
of $\ss_1$. That is, we partition the vector differently here.

When $(\mu_1, \sigma_1) = (\mu_2, \sigma_2)$ in Theorem \ref{theorem-2},
the asymptotic expansion of the likelihood ratio is an expansion for regular models:
\begin{equation}
2\{
\ell_n(\check \btheta^*) - \ell_n(\btheta_0) \}
=
n^{-1} (\sum_{i=1}^n A_i)^\tau \Sigma_A^{-1} (\sum_{i=1}^n A_i) +o_p(1).
\label{lik.dif.02}
\end{equation}
%
The result of Lemma \ref{lemma2} remains applicable:
\[
2\{ \ell_n(\btheta) - \ell_n(\btheta_0) \}
\leq
2 \ss_1^\tau \sum_{i=1}^n A_i + 2 \ss_2^\tau \sum_{i=1}^n B_i
-n\{ \ss_1^\tau \Sigma_A \ss_1  + \ss_2^\tau \Sigma_B \ss_2 \} + o_p(1).
\]
Since $\sigma_1= \sigma_2$
in Theorem \ref{theorem-2}(a), we have
\[
\ss_1
=(\mu, ~\sigma_+^2-1, \beta_0^2 + (\eta^2-1))^\tau; ~~
\ss_2  = (0, 0, \beta_0^4)^\tau.
\]
This leads to
\be
2\{ \ell_n(\tilde \btheta^*) - \ell_n(\btheta_0) \}
\leq
n^{-1} (\sum_{i=1}^n A_i)^\tau \Sigma_A^{-1} (\sum_{i=1}^n A_i)
+  6 n^{-1} \{ (\sum_{i=1}^n B_i[3])^{+}\}^2 +
o_p(1),
\label{lik.dif.21}
\ee
where we have $(\sum_{i=1}^n B_i[3])^{+}$
instead of $(\sum_{i=1}^n B_i[3])$ because of the
intrinsic constraint $\ss_2[3]= \beta^4_0 \geq 0$.
We skip the step of showing that the above upper bound
is attainable, since this is now routine.

Combining (\ref{lik.dif.02}) and (\ref{lik.dif.21}) gives
$$
R_{n,1}^*= 6 n^{-1} \{ (\sum_{i=1}^n B_i[3])^{+}\}^2 + o_p(1),
$$
which converges to $0.5\chi^2_0+0.5\chi^2_1$ in distribution,
which is conclusion (a).

For $R_{n,2}^*$ in (b), we are not helped by $\sigma_1 = \sigma_2$.
Yet
\[
2\{ \ell_n(\btheta) - \ell_n(\btheta_0) \}
\leq
2 \ss_1^\tau \sum_{i=1}^n A_i + 2 \ss_2^\tau \sum_{i=1}^n B_i
-n\{ \ss_1^\tau \Sigma_A \ss_1  + \ss^2 \Sigma_B \ss_2 \} + o_p(1)
\]
remains true for $\btheta$ in a small neighborhood of $\btheta_0$.
Similarly, we still have
\[
2\{ \ell_n(\hat \btheta^*) - \ell_n(\btheta_0) \}
=n^{-1} (\sum_{i=1}^n A_i)^\tau \Sigma_A^{-1} (\sum_{i=1}^n A_i)
+
\sup_{\ss_2}
\{ 2 \ss_2^\tau \sum_{i=1}^n B_i -n \ss_2^\tau \Sigma_B \ss_2  \}
+ o_p(1).
\]
We skip the proof that this upper bound is attained.
Hence,
\begin{equation}
\label{Rn2.star1}
R_{n,2}^*=\sup_{\ss_2}
\{ 2 \ss_2^\tau \sum_{i=1}^n B_i -n \ss_2^\tau \Sigma_B \ss_2  \}
+ o_p(1).
\end{equation}

The challenge is to provide an analytical description of the limiting
distribution when
\[
\ss_2  = (\beta_1^2, ~\beta_0\beta_1, ~\beta_0^4)^\tau.
\]
For this purpose, we
highlight the fact that $n^{-1/2}\sum_{i=1}^n B_i$ is asymptotically
multivariate normal with mean 0 and covariance matrix $\Sigma_B
= \mbox{diag}(1, 2, 1/6)$.
The supremum is hence attained in the range of $\ss_2 = O_p(n^{-1/2})$.
In the subregion where $|\beta_0| < n^{-1/7} = o(n^{-1/8})$,
we have $\ss_2[3] = \beta_0^8 < n^{-8/7} = o(n^{-1})$.
Hence,
\begin{eqnarray}
&& \hspace*{-5em}
\sup_{\ss_2, |\beta_0| < n^{-1/7}}
\{ 2 \ss_2^\tau \sum_{i=1}^n B_i -n \ss_2^\tau \Sigma_B \ss_2  \} \nonumber
\\
&=&
\sup_{\ss_2,  \beta_0=0}
\{ 2 \ss_2^\tau \sum_{i=1}^n B_i -n \ss_2^\tau \Sigma_B \ss_2  \}
+ o_p(1\nonumber)\\
&=&
n^{-1}\{ (\sum_{i=1}^n B_i[1])^+\}^2
+
(1/2) n^{-1}\{ \sum_{i=1}^n B_i[2]\}^2
+o_p(1).\label{Rn2.star2}
\end{eqnarray}
In the other subregion where $|\beta_0| \geq n^{-1/7}$,
combined with the restriction $\beta_0 \beta_1 = O_p(n^{-1/2})$,
we must have $\beta_1 = O_p(n^{-1/3})$.
Consequently, in this region,
$\ss_2[1] = \beta_1^2 = O(n^{-2/3})$.
This leads to
\[
\ss_2 [1] \sum_{i=1}^n B_i[1]  - n \{\ss_2 [1]\}^2 = o_p(1).
\]
Hence,
\begin{eqnarray}
&& \hspace*{-5em}
\sup_{\ss_2, |\beta_0| \geq n^{-1/7}}
\Big
\{ 2 \ss_2^\tau \sum_{i=1}^n B_i -n \ss_2^\tau \Sigma_B \ss_2
\Big\} \nonumber
\\
&=&
\sup_{\ss_2}
\Big \{
2\ss_2[2] \sum_{i=1}^n B_i[2]  - 2 n \ss_2[2]^2
+
2\ss_2[3] \sum_{i=1}^n B_i[3]  - (1/6) n \ss_2[3]^2
\Big  \}
+ o_p(1)\nonumber\\
&=&
(1/2)n^{-1}\{ \sum_{i=1}^n B_i[2]\}^2
+
6n^{-1}\{ (\sum_{i=1}^n B_i[3])^+\}^2
+o_p(1).
\label{Rn2.star3}
\end{eqnarray}

Combining (\ref{Rn2.star1})--(\ref{Rn2.star3}), we find
$$
R_{n,2}^*
=
(1/2)n^{-1}\{ \sum_{i=1}^n B_i[2]\}^2
+
\max \Big [
n^{-1}\{ (\sum_{i=1}^n B_i[1])^+\}^2,~
6n^{-1}\{ (\sum_{i=1}^n B_i[3])^+\}^2
\Big ]
+o_p(1).
$$
Therefore,  $R_{n,2}^*$
has the limiting distribution as claimed.

\section*{Acknowledgements}
The research is supported in part by
NSERC Grants RGPIN-2014-03743 and RGPIN-2015-06592 
and Singapore Ministry Education Academic Research Fund Tier 1 
and the Ministry of Education of Singapore: MOE2014-T2-1- 072.

%


\end{document}